\documentclass[12pt]{book}

\makeatletter

\newdimen\earraycolsep
\def\eqnarray{\stepcounter{equation}\let\@currentlabel\theequation
\global\@eqnswtrue\m@th
\global\@eqcnt\z@\tabskip\@centering\let \\\@eqncr
$$\halign to\displaywidth\bgroup\@eqnsel\hskip\@centering
$\displaystyle\tabskip\z@{##}$&\global\@eqcnt\@ne
\hskip 2\earraycolsep \hfil$\displaystyle{##}$\hfil
&\global\@eqcnt\tw@ \hskip 2\earraycolsep $\displaystyle\tabskip\z@
{##}$\hfil
\tabskip\@centering&\llap{##}\tabskip\z@\cr}
\renewcommand{\theequation}{\arabic{equation}}
\renewcommand{\thetable}{\arabic{table}}
\renewcommand{\thefigure}{\arabic{figure}}

\def\title{\chapter}
\renewcommand\chapter{\ifodd\c@page\clearpage\else
          \cleardoublepage\fi
          \global\@topnum\z@
          \@afterindenttrue
          \secdef\@chapter\@schapter}
\def\@makechapterhead#1{%
  \vspace*{120\p@}%
  {\parindent \z@ \raggedright \reset@font
    \bfseries #1\par
    \nobreak
    \vskip 36\p@
  }}
\def\author#1{{\pretolerance=10000 \raggedright \advance \leftskip
by 1in \noindent #1 \vskip 1pc}}
\def\affiliation#1{{\advance\leftskip by 1in \noindent #1 \vskip
-1pc}}
\def\refnote#1{{$^{\hbox{\scriptsize #1}}$}}

\renewcommand\section{\@startsection{section}{1}{\z@}{2pc \@plus
1ex minus
    .2ex}{1pc \@plus .2ex}{\reset@font\normalsize\bfseries}}
\renewcommand\subsection{\@startsection{subsection}{2}{\z@}{1pc
\@plus 1ex
    minus.2ex}{1pc \@plus .2ex}{\reset@font\normalsize\bfseries}}
\renewcommand\subsubsection{\@startsection{subsubsection}{3}
{\parindent}
   {1pc \@plus 1ex minus.2ex}{-0.5em}{\reset@font\normalsize
\bfseries}}

\def\AmS{{\protect\the\textfont2 A\kern-.1667em\lower.5ex\hbox
{M}\kern-.125emS}}

\def\p@LaTeX{{\family{times}\series{m}\shape{n}\selectfont L\kern
-.36em\raise.3ex\hbox{\scriptsize A}\kern-.15em
T\kern-.1667em\lower.7ex\hbox{E}\kern-.125emX}}

\newlength{\colwidth}

\def\@oddhead{\hfil}
\def\@evenhead{\hfil}
\def\@oddfoot{{\bfseries\hfil\thepage}}
\def\@evenfoot{{\bfseries\thepage\hfil}}
\def\fnum@figure{\footnotesize\raggedright{\bfseries \figurename
~\thefigure.}}
\def\fnum@table{\normalsize\raggedright{\bfseries \tablename
~\thetable.}}
\long\def\@makecaption#1#2{\vskip 10\p@ {#1 #2\par}}
\long\def\@makefntext#1{\setbox0=\hbox{$\m@th^{\@thefnmark}$}
\noindent\hangindent=\wd0 \box0 #1}
\newbox\@atbox
\long\def\atable#1#2#3{\begin{table}[tbp]\centering\footnotesize
\setbox\@atbox\hbox{#2}
\parbox{\wd\@atbox}{\caption{#1}}\par\smallskip
#2
\par\smallskip\parbox{\wd\@atbox}{\raggedright #3}
\end{table}}
\def\@bibitem{\par\noindent \hangindent=2pc \hangafter=1}
\def\thebibliography{%
\section*{REFERENCES}%
\bgroup\footnotesize
\def\newblock{\hskip .11em plus.33em minus.07em}%
\let\bibitem\@bibitem}
\def\endthebibliography{\par\egroup}
\def\@nbibitem#1{\par\noindent \hangindent=2pc \hangafter=1
\refstepcounter{enumi}\hbox to 2pc{\arabic{enumi}.\hfil}%
\immediate\write\@auxout{\string\bibcite{#1}{\arabic{enumi}}}}
\def\numbibliography{%
\section*{REFERENCES}%
\bgroup\footnotesize
\setcounter{enumi}{0}%
\def\newblock{\hskip .11em plus.33em minus.07em}%
\let\bibitem\@nbibitem}
\def\endnumbibliography{\par\egroup}
\def\@cite#1#2{{#1\if@tempswa , #2\fi}}

\makeatother

\def\req#1{(\ref{#1})}

\def\setfonts{%
\font\frbig=eufm10 scaled\magstep1
\font\frscr=eufm8 
\font\frscrscr=eufm8
\newfam\frfam
\textfont\frfam=\frbig
\scriptfont\frfam=\frscr
\scriptscriptfont\frfam=\frscrscr
\def\fr{\fam\frfam}

\font\openbig=msbm10 scaled\magstep1
\font\openscr=msbm8 scaled\magstephalf
\font\openscrscr=msbm8
\newfam\openfam
\textfont\openfam=\openbig
\scriptfont\openfam=\openscr
\scriptscriptfont\openfam=\openscrscr
\def\open{\fam\openfam}

\font\ssfbig=cmss10 scaled\magstep1
\font\ssfscr=cmss8 scaled\magstephalf
\font\ssfscrscr=cmss8
\newfam\ssffam
\textfont\ssffam=\ssfbig
\scriptfont\ssffam=\ssfscr
\scriptscriptfont\ssffam=\ssfscrscr
\def\ssf{\fam\ssffam}
}

\makeatletter
\@addtoreset{equation}{section}
\newdimen\normalarrayskip
\newdimen\minarrayskip
\normalarrayskip\baselineskip
\minarrayskip\jot
\newif\ifold \oldtrue \def\new{\oldfalse}
\def\arraymode{\ifold\relax\else\displaystyle\fi}

\def\@arrayskip{\ifold\baselineskip\z@\lineskip\z@
  \else
  \baselineskip\minarrayskip\lineskip2\minarrayskip\fi}
\def\@arrayclassz{\ifcase \@lastchclass \@acolampacol \or
\@ampacol \or \or \or \@addamp \or
 \@acolampacol \or \@firstampfalse \@acol \fi
\edef\@preamble{\@preamble
 \ifcase \@chnum
  \hfil$\relax\arraymode\@sharp$\hfil
  \or $\relax\arraymode\@sharp$\hfil
  \or \hfil$\relax\arraymode\@sharp$\fi}}
\def\@array[#1]#2{\setbox\@arstrutbox=\hbox{\vrule
  height\arraystretch \ht\strutbox
  depth\arraystretch \dp\strutbox
  width\z@}\@mkpream{#2}\edef\@preamble{\halign \noexpand\@halignto
\bgroup \tabskip\z@ \@arstrut \@preamble \tabskip\z@ \cr}%
\let\@startpbox\@@startpbox \let\@endpbox\@@endpbox
 \if #1t\vtop \else \if#1b\vbox \else \vcenter \fi\fi
 \bgroup \let\par\relax
 \let\@sharp##\let\protect\relax
 \@arrayskip\@preamble}
\@addtoreset{equation}{section}
\makeatother

\setfonts

\def\d{\partial}
\def\spsi{\psi^*}

\def\N#1{N\!=\!#1}
\def\tSl#1{{\widehat{s\ell}}(#1)}
\def\tSL#1{{\widehat{s\ell}}(#1)}
\def\SL#1{s\ell(#1)}

\def\half{\frac{1}{2}}

\def\cG{{\cal G}}
\def\cH{{\cal H}}

\def\cL{{\cal L}}

\def\cO{{\cal O}}
\def\cP{{\cal P}}
\def\cQ{{\cal Q}}

\def\cT{{\cal T}}
\def\cU{{\cal U}}

\def\oN{{\open N}}
\def\oC{{\open C}}
\def\oQ{{\open Q}}

\def\oZ{{\open Z}}

\def\ctop{{\ssf c}}
\def\Ctop{{\ssf C}}
\def\htop{{\ssf h}}

\def\hplus{{\ssf h}^+}
\def\hminus{{\ssf h}^-}

\def\Lambdach{\Lambda_{\rm ch}}

\def\tensor{\otimes}

\def\mm{\cal}
\def\smm{\fr}

\def\mM{{\mm M}}

\def\mF{{\mm F}}

\def\mR{{\mm R}}

\def\mU{{\mm U}}
\def\mV{{\mm V}}

\def\smM{{\smm M}}

\def\smV{{\smm V}}
\def\smR{{\smm R}}

\def\bar{\overline}
\def\frac#1#2{\mathchoice{{\textstyle{{#1}\over{#2}}}}{{#1\over#2}}%
  {{#1\over#2}}{{#1\over#2}}}
\def\ket#1{\mathchoice{%
{\left|{#1}\right\rangle}}{|{#1}\rangle}{|{#1}\rangle}{|{#1}\rangle}}
\def\kettop#1{\mathchoice{{\left|{#1}\right\rangle}_{\rm top}}%
  {|{#1}\rangle_{\rm top}}{|{#1}\rangle_{\rm top}}{|{#1}\rangle_{\rm top}}}
\def\ketSL#1{\left|{#1}\right\rangle_{\SL2}}

\def\ketstar#1{\mathchoice{%
    {\left|{#1}\right\rangle^{\phantom{y}}_{*}}}{|{#1}\rangle_{*}}{|{#1}\rangle_{*}}{|{#1}\rangle_{*}}}

\def\ketm#1{\mathchoice{%
    {\left|{#1}\right\rangle_{\rm m}}}{|{#1}\rangle_{\rm m}}{|{#1}\rangle_{\rm m}}{|{#1}\rangle_{\rm m}}}

\def\ketGH#1{\left|{#1}\right\rangle_{\rm GH}}

\newtheorem{thm}{Theorem}

\def\VER{{\cal VER}}

\def\TOP{{\cal TOP}}

\def\emt{energy-momentum tensor}
\def\hw{highest-weight}

\def\NPB{Nucl.\ Phys.\ B}
\def\PLB{Phys.\ Lett.\ B}
\def\MPLA{Mod.\ Phys.\ Lett.\ A}

\def\IJMPA{Int.\ J.\ Mod.\ Phys.\ A}

\begin{document}

\begin{flushright}
  {\tt q-alg/9712024}
\end{flushright}

\vspace*{40pt}

\noindent
{{\bfseries PAST THE HIGHEST-WEIGHT, AND WHAT YOU CAN FIND
    THERE}\footnote{Talk presented at the NATO Advanced Research
    Workshop on Theoretical Physics `New Developments in Quantum Field
    Theory', June 14--20, 1997, Zakopane, and IV International
    Conference `Conformal Field Theories and Integrable Models',
    Chernogolovka, June 23--27, 1997.}}\nobreak\vskip 36pt

{\advance \leftskip by 1in \noindent \large A.M.~Semikhatov \vskip 1pc}

{\advance \leftskip by 1in \noindent {\sl Lebedev Physics Institute,
    Moscow 117924, Russia}\vskip-1pc}

\vskip30pt

\parbox{.9\textwidth}{\footnotesize The properties of highest-weight
  representations of the $\N2$ superconformal algebra in two
  dimensions can be considerably simplified when re-expressed in terms
  of relaxed $\tSL2$ representations. This applies to the appearance
  of submodules and hence, of singular vectors, and to the structure
  of the embedding diagrams and the BGG-type resolution. I also discuss
  the realization of these representations in the bosonic string,
  where the generalized DDK prescription amounts to the requirement
  that the representations have a charged singular vector, and the
  role of the fermionic screening operator.}

\vspace*{-4pt}

\section{INTRODUCTION}

\vspace*{-4pt}

In this talk, I discuss how the crucial features of modules over the
$\N2$ superconformal algebra in two dimensions reformulate in a
simpler way in terms of modules over the affine algebra~$\tSl2$. The
key statement, which can be advertised as\strut

\centerline{\framebox{\it the $\N2$ and affine-$\SL2$ representations
    theories are ``essentially equivalent''}}

\noindent
was proved in\,\refnote{\cite{[FST]}}.\strut \ Also analysed in that
paper were the degenerations (reducibility patterns) of the $\tSL2$
representations (called the {\sc relaxed} Verma modules) that `model'
the $\N2$ Verma-like modules. In a parallel development, degenerations
of the $\N2$ Verma modules were directly analysed
in\,\refnote{\cite{[ST4]}} (and, as a by-product, the problem of $\N2$
{\it sub\/}singular vectors was resolved).

In this talk, I will be much less formal than in those papers, and I
will also discuss briefly some new
developments\,\refnote{\cite{[SSi],[q]}} stemming from the results
of~\refnote{\cite{[FST],[ST4]}}. These are aimed at constructing the
BGG-type resolution\,\refnote{\cite{[BGG]}} for the irreducible
representations and, thus, at systematically deriving the characters
and, on the other hand, at deriving the complete quantum symmetry of
the $\tSL2$ fusion rules\,\refnote{\cite{[AY],[FM],[Andreev],[PRY]}}.
Before proceeding to more precise formulations, let us see whether the
``essential equivalence'' of the $\N2$ and $\tSL2$ representation
theories comes as a news: 
\begin{list}{}{\leftmargin12pt}\addtolength{\parskip}{-9pt}

\item[--] on the one hand, the two algebras appear to have little in
  common, since one is a rank-2 (bosonic) affine Lie algebra, while
  the other is a rank-3 {\it super\/}algebra that is {\it not\/} an
  affine Lie algebra;

\item[--] on the other hand, the appearance of the two algebras in CFT
  is often `correlated', they `share' parafermionic theories, etc.
  \pagebreak[3]
\end{list}

\vspace{-4pt}

\noindent
Surprisingly or not, establishing the equivalence requires introducing
somewhat unusual $\tSL2$ representations; more precisely,
\begin{list}{}{\leftmargin16pt}\addtolength{\parskip}{-9pt}
\item[1.] One considers an arbitrary complex level
  $k\in\oC\setminus\{-2\}$ on the $\tSL2$ side.

\item[2.] On the $\N2$ side, one considers the `standard'
  representation category, which includes the Verma modules and their
  quotients, the (unitary and non-unitary) irreducible
  representations,~etc., along with their images under the spectral
  flow ({\it twists\/}), for the central
  charge~$\ctop\in\oC\setminus\{3\}$.

\item[3.] {\it Modulo the spectral flow transform}, this $\N2$
  category is equivalent to the category of $\tSL2$ representations
  {\it of the type that has not been considered before\/} --- the {\sc
    relaxed} \hw-type representations (and their twists).

\item[4.] On the other hand, the standard \hw-type $\tSL2$
  representations (category $\cO$\,\refnote{\cite{[TheBook]}} {\it and
    their twists\/}) turn out to be related to a narrower category of
  (twisted) {\sc topological} $\N2$ Verma modules.
\end{list}

\vspace*{-4pt}

A statement regarding the equivalence of two categories can often be
interpreted to the effect that there are two different languages to
describe the same structure. In this way, any `structural' result that
is claimed about $\N2$ Verma modules can in principle be seen in
relaxed $\tSL2$ modules, and vice versa.  As it may (and does) happen
with equivalence of categories, however, an exciting point is that a
number of fairly obvious facts about $\tSL2$ representations translate
into the statements which are not quite obvious for the $\N2$
representations.  Thus, we see that the objects that are quite
standard on the $\N2$ side can be described in the $\tSL2$ terms by
introducing a new type of modules, while only a subclass of $\N2$
representations corresponds to the standard $\tSL2$ representations.
Of a crucial importance is, therefore, the distinction between two
different types of Verma-like modules over each algebra; for the $\N2$
algebra, this distinction is masked due to an effect that we are going
to discuss, while from the $\tSL2$ point of view, the distinction is
much easier to see, and it can roughly be summarized by saying that in
the relaxed Verma modules,\strut

\centerline{\framebox{\it one goes past the highest-weight vector,}}

\noindent
which are going to describe in more detail.\strut \ As a good
illustration to the equivalence of categories, let me note that, while
the $\tSL2$ representations where `one goes past the \hw{}' may look
somewhat unusual, {\sc the same\/} effects described in the $\N2$
context are not considered unusual at all!

\section{$\tSL2$ HIGHEST-WEIGHT REPRESENTATIONS}

Let us begin with the $\tSL2$ algebra. We fix the level $k\neq-2$.
Recall what one does when constructing a \hw-type module. The
generators are broken into, roughly, two `halves', one of which are
declared annihilation operators with respect to a {\it \hw{} vector},
while the others {\it create\/} states, except for the `Cartan'
generator(s), whose eigenvalues simply `label' the \hw{} vectors:
\begin{equation}
  J^+_{\geq0}\,\ketSL{j,k}=J^0_{\geq1}\,\ketSL{j,k}=
  J^-_{\geq1}\,\ketSL{j,k}=0\,,\quad
  J^0_{0}\,\ketSL{j,k}=j\,\ketSL{j,k}\,,
  \label{sl2hig}
\end{equation}
where $j,\,k\in\oC$.  In the {\it Verma\/} module, by definition,
there are no relations among the states produced by the creation
operators from the \hw{} vector.  That is, the Verma module is {\it
  freely\/} generated from the \hw{} vector by the operators declared
to be the creation ones.  The structure of $\tSL2$ Verma modules is
conveniently encoded in the {\it extremal diagram\/}
\begin{equation}
  \unitlength=1pt
  \begin{picture}(250,60)
    \put(0,-15){
    \put(-35,62){\Huge $\ldots$}
    \put(0,60){$\bullet$}
    \put(15,65){${}^{J^-_0}$}
    \put(28,63){\vector(-1,0){22}}
    \put(30,60){$\bullet$}
    \put(45,65){${}^{J^-_0}$}
    \put(58,63){\vector(-1,0){22}}
    \put(60,60){$\bullet$}
    \put(75,65){${}^{J^-_0}$}
    \put(88,63){\vector(-1,0){22}}
    \put(90,60){$\circ$}
    \put(103,64){${}_{J^+_{-1}}$}
    \put(97,60){\vector(2,-1){17}}
    \put(115,47){$\bullet$}
    \put(24,-13){%
      \put(103,64){${}_{J^+_{-1}}$}
      \put(97,60){\vector(2,-1){17}}
      \put(115,47){$\bullet$}
      }
    \put(48,-26){%
      \put(103,64){${}_{J^+_{-1}}$}
      \put(97,60){\vector(2,-1){17}}
      \put(115,47){$\bullet$}
      }
    \put(72,-39){%
      \put(103,50){\LARGE $\cdot$}
      \put(109,47){\LARGE $\cdot$}
      \put(115,44){\LARGE $\cdot$}
      }
    }
  \end{picture}
  \label{Vermaextr}
\end{equation}
The states shown in the diagram are {\it extremal\/} in the sense that
they have boundary values of the $({\rm charge}, {\rm level})$
bigrading; all of the other states of the module should be thought of
as lying in the interior of the angle in the diagram (on the
rectangular lattice according to their (charge, level)). Finding a
{\it submodule\/} in the Verma module can be (somewhat more
schematically) represented as
\begin{equation}
  \unitlength=1pt
  \begin{picture}(260,36)
    \put(0,-15){
      \put(10,50){\line(1,0){120}}
      \put(130,50){\line(2,-1){60}}
      \put(15,33){\line(1,0){90}}
      \put(105,33){\line(2,-1){30}}
      }
  \end{picture}
\end{equation}
Whenever one considers quotients of Verma modules, the extremal
diagrams become `smaller', as some of the states are eliminated from
the module. All of such extremal diagrams, therefore, satisfy the
following criterion:\strut\\
\centerline{\kern-8pt\framebox{\kern-4pt\it Any straight line going
    through any state intersects the boundary on at least
    one end,}\kern-4pt}\\
which is formalized as follows:\strut\ for any state $\ket{X}$ from the
module,
\begin{equation}
  \forall n\in\oZ,\quad\exists N\in\oN: \quad
  ({\rm either}\quad (J^+_{n})^{N}\,\ket X=0 \quad{\rm or}\quad
  (J^-_{-n})^{N}\,\ket X=0)
  \label{Jterminate}
\end{equation}
However, this criterion also selects the so-called {\it twisted\/}
modules. As regards twisted Verma modules, their extremal diagrams are
`rotations' of the above, e.g.:
\begin{equation}
  \unitlength=1pt
  \begin{picture}(250,60)
    \put(0,0){\line(1,2){20}}
    \put(20,40){\line(2,1){40}}
    \put(100,30){or}
    \put(200,50){\line(1,0){80}}
    \put(140,20){\line(2,1){60}}
  \end{picture}
\end{equation}
In more formal terms, a {\it twisted Verma module\/}
$\smM_{j,k;\theta}$ is freely generated by $J^+_{\leq\theta-1}$,
$J^-_{\leq-\theta}$, and $J^0_{\leq-1}$ from a {\it twisted
  highest-weight vector\/} $\ketSL{j,k;\theta}$ defined by the
conditions
\begin{equation}\new
  \begin{array}{l}
    J^+_{\geq\theta}\,\ketSL{j,k;\theta}=J^0_{\geq1}\,\ketSL{j,k;\theta}=
    J^-_{\geq-\theta+1}\,\ketSL{j,k;\theta}=0\,,\\
    \left(J^0_{0}+\frac{k}{2}\theta\right)\,\ketSL{j,k;\theta}=
    j\,\ketSL{j,k;\theta}.
  \end{array}
  \label{sl2higgeneral}
\end{equation}
The mapping that underlies the construction of twisted modules is
known as the spectral flow transform\refnote{\cite{[BH]}}
\begin{equation}
  \cU_\theta:\quad
  J^+_n\mapsto J^+_{n+\theta}\,,\qquad
  J^-_n\mapsto J^-_{n-\theta}\,,\qquad
  J^0_n\mapsto J^0_n+\frac{k}{2}\theta\delta_{n,0}\,,\qquad\theta\in\oZ\,.
  \label{spectralsl2}
\end{equation}
Then,
\begin{equation}
  \cU_\theta\,\ketSL{j,k}=\ketSL{j,k;\theta}\,,\qquad
  \ketSL{j,k}=\ketSL{j,k;0}\,.
\end{equation}

A characteristic feature of extremal diagrams of $\tSL2$ Verma modules
is the existence of the `angle' that corresponds to the \hw{} vector.
This will reappear in the topological $\N2$ Verma modules. We now
consider the second main ingredient, the $\N2$ algebra and
representations, after which we return to $\tSL2$ and ``go past the
\hw''; those representations won't have the angle in the extremal
diagrams.

\section{$\N2$ ALGEBRA AND REPRESENTATIONS}

The $\N2$ superconformal algebra contains two fermionic currents,
$\cQ$ and $\cG$, in addition to the Virasoro generators $\cL$ and the
$U(1)$ current $\cH$. The commutation relations read as
\begin{equation}\new
  \begin{array}{lclclcl}
    \left[\cL_m,\cL_n\right]&=&(m-n)\cL_{m+n}\,,&\qquad&[\cH_m,\cH_n]&=
    &\frac{\Ctop}{3}m\delta_{m+n,0}\,,\\

    \left[\cL_m,\cG_n\right]&=&(m-n)\cG_{m+n}\,,&\qquad&[\cH_m,\cG_n]&=&\cG_{m+n}\,,
    \\
    \left[\cL_m,\cQ_n\right]&=&-n\cQ_{m+n}\,,&\qquad&[\cH_m,\cQ_n]&=&-\cQ_{m+n}\,,\\

    \left[\cL_m,\cH_n\right]&=&\multicolumn{5}{l}{-n\cH_{m+n}+\frac{\Ctop}{6}(m^2+m)
      \delta_{m+n,0}\,,}\\
    \left\{\cG_m,\cQ_n\right\}&=&\multicolumn{5}{l}{2\cL_{m+n}-2n\cH_{m+n}+
      \frac{\Ctop}{3}(m^2+m)\delta_{m+n,0}\,,}
  \end{array}\qquad m,~n\in\oZ\,.
  \label{topalgebra}
\end{equation}
The element $\Ctop$ is central; in representations, we will not
distinguish between $\Ctop$ and its eigenvalue $\ctop\!\in\!\oC$,
which it will be convenient to parametrize as $\ctop=3\,\frac{t-2}{t}$
with $t\!\in\!\oC\setminus\{0\}$. The $\N2$ spectral flow
transform\,\refnote{\cite{[SS]}} acts as follows:
\begin{equation}
  {\cal U}_\theta:\new
  \begin{array}{rclcrcl}
    \cL_n&\mapsto&\cL_n+\theta\cH_n+\frac{\ctop}{6}(\theta^2+\theta)
    \delta_{n,0}\,,&{}&
    \cH_n&\mapsto&\cH_n+\frac{\ctop}{3}\theta\delta_{n,0}\,,\\
    \cQ_n&\mapsto&\cQ_{n-\theta}\,,&{}&\cG_n&\mapsto&\cG_{n+\theta}\,,
  \end{array}
  \label{U}
\end{equation}

\subsection{Massive and topological $\N2$ modules}


There exist two types of Verma-like modules over the $\N2$ algebra,
which we call the massive and the topological ones; in the literature,
the former are commonly considered as `the' $\N2$ Verma modules, while
the latter are precisely those $\N2$ modules that are in a good
correspondence with the $\tSL2$ Verma modules from the previous
section.

A {\it massive Verma module\/} $\mU_{h,\ell,t}$ is freely generated by
the generators $\cL_{-m}$, $\cH_{-m}$, $\cG_{-m}$, $m\in\oN$, and
$\cQ_{-m}$, $m\in\oN_0$ (with $\oN\!=\!1,2,\ldots$ and $\oN_0\!=\!0,
1,2,\ldots$) from a {\it massive \hw{} vector\/} $\ket{h,\ell,t}$
satisfying the following set of highest-weight conditions:
\begin{equation}\new
  \begin{array}{rcl}
    \cQ_{\geq1}\,\ket{h,\ell,t}\kern-4pt&=&
    \kern-4pt\cG_{\geq0}\,
    \ket{h,\ell,t}= \cL_{\geq 1}\,\ket{h,\ell,t}=
    \cH_{\geq1}\,\ket{h,\ell,t}=0\,,\\
    \cH_0\,\ket{h,\ell,t}\kern-4pt&=&
    \kern-4pt
    h\,\ket{h,\ell,t}\,,\qquad
    \cL_0\,\ket{h,\ell,t}
    = \ell\,\ket{h,\ell,t}\,.
  \end{array}
  \label{masshw}
\end{equation}
In the bigrading implied by (charge,\,level), or more precisely,
by the eigenvalues of $(-\cH_0,\,\cL_0)$), the extremal diagram of a
massive Verma module has the shape of a parabola for the simple reason
that, having acted on the \hw{} vector with, say, $\cQ_0$, applying
the same operator once again gives identical zero, and `the best one
can do' to construct a state with the extremal (charge, level)
bigrading is to act with the $\cQ_{-1}$ mode, etc.:
\begin{equation}
  \unitlength=1.00mm
  \begin{picture}(140,40)
    \put(50.00,5.00){
      \put(00.00,00.00){$\bullet$}
      \put(10.00,20.00){$\bullet$}
      \put(10.00,20.00){$\bullet$}
      \put(20.00,30.00){$\bullet$}
      \put(30.00,30.00){$\bullet$}
      \put(40.00,20.00){$\bullet$}
      \put(50.00,00.00){$\bullet$}
      \put(10.00,18.70){\vector(-1,-2){8}}
      \put(19.50,29.50){\vector(-1,-1){7}}
      \put(22.40,31.10){\vector(1,0){7}}
      \put(32.50,29.80){\vector(1,-1){7}}
      \put(42.00,19.00){\vector(1,-2){8}}
      \put(00.00,13.00){${}_{\cG_{-2}}$}
      \put(09.00,26.50){${}_{\cG_{-1}}$}
      \put(23.90,33.50){${}_{\cQ_{0}}$}
      \put(37.00,27.50){${}_{\cQ_{-1}}$}
      \put(47.00,13.00){${}_{\cQ_{-2}}$}
      \put(11.00,32.00){${}_{\ket{h_,\ell,t}}$}
      \put(00.50,-06.00){$\vdots$}
      \put(50.50,-06.00){$\vdots$}
      }
  \end{picture}
  \label{massdiagramdouble}
\end{equation}
An important fact is that all of the states on the extremal diagram
satisfy the annihilation conditions
\begin{equation}
  \cQ_{-\theta+m+1}\approx \cG_{\theta+m}\approx
  \cL_{m+1}\approx\cH_{m+1}\approx0\,,\quad m\in\oN_0
  \label{twistedannihil}
\end{equation}
for $\theta$ ranging over the integers, from $-\infty$ in the left end
to $+\infty$ in the right end of the parabola.

Now, 
there can be {\it two\/} different types of Verma submodules
in~$\mU_{h,\ell,t}$. In the language of extremal diagrams, these look
like (with the discrete parabolas replaced by smooth ones for
simplicity)
\begin{equation}
  \unitlength=1pt
  \begin{picture}(350,80)
    \bezier{600}(0,0)(60,160)(120,0)
    \bezier{600}(20,0)(70,110)(110,0)
    \put(180,40){or}
    \put(230,0){
      \bezier{600}(0,0)(60,160)(120,0)
      \bezier{400}(10,25)(60,100)(100,0)
      \put(10,25){$\bullet$}
      }
  \end{picture}
  \label{pic:compare}
\end{equation}
In the first case, we have a {\it massive\/} Verma {\it sub\/}module,
all of the states on its extremal diagram (as well as those on the
extremal diagram of the module itself) satisfying the annihilation
conditions from~\req{masshw} for $\theta\in\oZ$.  In the other
case, on the contrary, there is a distinguished state, marked with a
$\bullet$, that satisfies the annihilation conditions
\begin{equation}
  \cQ_{-\theta+m}\approx \cG_{\theta+m}\approx
  \cL_{m+1}\approx\cH_{m+1}\approx0\,,\quad m\in\oN_0
  \label{annihiltop}
\end{equation}
for a fixed $\theta\in\oZ$. Such states will be referred to as
(twisted) {\it topological \hw{} vectors}, and in the above context,
as {\it topological singular vectors} (the $\theta=0$ case being the
`untwisted' one). The precise definition is as follows.

A {\it twisted topological Verma module\/} $\smV_{h,t;\theta}$ is
freely generated by $\cL_{\leq-1}$, $\cH_{\leq-1}$,
$\cG_{\leq\theta-1}$, and $\cQ_{\leq-\theta-1}$ from a twisted {\it
  topological \hw{} vector\/} subjected to annihilation
conditions~\req{annihiltop}, where, in addition,
\begin{equation}\new
  \begin{array}{rcl}
    (\cH_0+\frac{\ctop}{3}\theta)\,\ket{h,t;\theta}_{\rm
      top}&=&
    h\,\ket{h,t;\theta}_{\rm top}\,,\\
    (\cL_0+\theta\cH_0+\frac{\ctop}{6}(\theta^2+\theta))
    \,\ket{h,t;\theta}_{\rm top}&=&0\,.
  \end{array}
  \label{Cartantheta}
\end{equation}

A characteristic feature of the extremal diagram of a topological
Verma module is the existence of a `cusp', i.e.\ a state that
satisfies stronger (the twisted topological) \hw{} conditions than the
other states in the diagram. As a result, the extremal diagram is
narrower than that of a massive Verma module.  In the $\theta=0$ case
for simplicity, the extremal diagram of $\mV_{h,t}\equiv\smV_{h,t;0}$
reads as (with $\ket{h,t}_{\rm top}\equiv \ket{h,t;0}_{\rm top}$)
\begin{equation} \unitlength=0.9mm
  \begin{picture}(140,40)
    \put(50.00,5.0){
      \put(00.00,00.00){$\bullet$}
      \put(10.00,20.00){$\bullet$}
      \put(10.00,20.00){$\bullet$}
      \put(20.00,30.00){$\bullet$}
      \put(29.70,20.00){$\bullet$}
      \put(40.00,00.00){$\bullet$}
      \put(9.70,19.00){\vector(-1,-2){8}}
      \put(19.70,29.70){\vector(-1,-1){7}}
      \put(22.00,29.70){\vector(1,-1){7}}
      \put(32.00,19.00){\vector(1,-2){8}}
      \put(00.00,13.00){${}_{\cG_{-2}}$}
      \put(11.00,28.00){${}_{\cG_{-1}}$}
      \put(27.00,28.00){${}_{\cQ_{-1}}$}
      \put(37.00,13.00){${}_{\cQ_{-2}}$}
      \put(19.00,34.00){${}_{\ket{h,t}_{\rm top}}$}
      \put(00.50,-06.00){$\vdots$}
      \put(40.50,-06.00){$\vdots$}
      }
  \end{picture}
  \label{topdiag}
\end{equation}

When taking quotients, the extremal diagrams may only become smaller,
which allows us to formulate a criterion that automatically singles
out the {\it topological \hw{}-type\/} modules (the corresponding
$\cO$-category).  For any $n\in\oZ$, by the `massive' parabola
$\cP(n,X)$ running through a state $\ket{X}$, we understand the set of
states
\begin{equation}
  \cQ_{n-N}\,\ldots\,\cQ_{n-1}\,\cQ_{n}\,\ket X=0\,,\quad
  \cG_{-n-M}\,\ldots\,\cG_{-n-2}\,\cG_{-n-1}\,\ket X=0\,,
  \qquad N, M\in\oN\,.
  \label{terminate}
\end{equation}
Then, a module belongs to the twisted `topological'
$\cO$-category if, for any state $\ket{X}$,\strut\\
\centerline{\kern-3pt\framebox{\kern-3pt\it any massive parabola
    intersects
    the extremal diagram of the module on at least one end}}\\
which, again, means simply that the states \req{terminate} become zero
in at least one branch, either for $N\gg1$ or for $M\gg1$.

The massive $\N2$ Verma modules do not satisfy this criterion.
However, in the massive case as well, one can formulate a criterion
that does not allow the modules to become too wide: for any element
$\ket{X}$,
\begin{equation}
  \forall n\in\oZ,\quad
    {\rm either}~\ldots\cQ_{n-2}\,\cQ_{n-1}\,\cQ_{n}\,\ket X=0\quad
    {\rm or}~\ldots\cG_{-n-2}\,\cG_{-n-1}\,\cG_{-n}\,\ket X=0\,.
\end{equation}

In the next Section, we address the problem of finding the $\tSL2$
counterpart of the above $\N2$ modules.  We first map the generators
and then investigate the representations.

\vspace{-6pt}

\section{FROM $\N2$ TO $\tSL2$}

\vspace{-6pt}

\subsection{An operator construction}

\vspace{-6pt}

We now use an operator construction allowing us to build the $\tSL2$
currents out of the $\N2$ generators and a free scalar with the
operator product $\phi(z)\phi(w)=-\ln(z-w)$.  As a necessary
preparation, we `pack' the modes of the $\N2$ generators into the
corresponding fields, $\cT(z)=\sum_{n\in\oZ}\cL_n z^{-n-2}$, \
$\cG(z)=\sum_{n\in\oZ}\cG_n z^{-n-2}$, $\cQ(z)=\sum_{n\in\oZ}\cQ_n
z^{-n-1}$, and $\cH(z)=\sum_{n\in\oZ}\cH_n z^{-n-1}$, and similarly
with the $\tSL2$ currents. We also define vertex operators
$\psi=e^\phi$ and $\spsi=e^{-\phi}$.  Then, for $\ctop\neq3$,
\begin{equation}
    J^+= \cQ\psi\,,\qquad J^-=\frac{3}{3-\ctop}\, \cG\spsi\,,\qquad
    J^0=-\frac{3}{3-\ctop}\,\cH+\frac{\ctop}{3-\ctop}\,\d\phi
  \label{invKS}
\end{equation}
are the $\tSL2$ generators of level
$k={2\ctop\over3-\ctop}$.\footnote{At the Conference, M.~Halpern told
  me that such a mapping had been known to M.~Peskin {\it et al.}, but
  I could not find the reference.} \

We also have a free scalar with
signature $-1$, whose modes commute with the $\tSL2$
generators~\req{invKS}:
\begin{equation}
  I^-=\sqrt{(k+2)/2}\,(\cH-\d\phi)\,.
  \label{dF}
\end{equation}
The modes $I^-_n$ generate a Heisenberg algebra. Then the module
$\mF^-_q$ is defined as a Verma module over this Heisenberg algebra
with the \hw{} vector defined by \ $I^-_n\ket{q}^-=0$, $n\geq1$, and
$I^-_0\ket{q}^-=q\ket{q}^-$.

\subsection{Relating the representations}

The behaviour of {\it representations\/} under operator constructions
of this sort can be quite complicated.\footnote{Recall, for example,
  how the $\tSL2$ Verma modules are rearranged under the Wakimoto
  bosonization\refnote{\cite{[FFr]}} --- Wakimoto modules more or less
  `interpolate' between Verma and contragredient Verma modules.} \ In
our case, we take a topological Verma module $\mV_{h,t}$ and {\it
  tensor\/} it with the module~$\Xi$ of the free scalar.  The latter
module is defined as $\Xi=\oplus_{n\in\oZ}\mF_n$, where $\mF_n$ is a
Verma module with the highest-weight vector $\ket{n}_\phi$ such that
\begin{equation}
    \phi_m\ket{n}_\phi=0\,,~ m\geq1\,,\quad
    \psi_m\ket{n}_\phi=0\,,~ m\geq n+1\,,\quad
    \spsi_m\ket{n}_\phi=0\,,~ m\geq -n+2\,,
\end{equation}
and $\phi_0\ket{n}_\phi=-n\ket{n}_\phi$. We then have the following
Theorem:

\vspace*{-4pt}

\begin{thm}[\refnote{\cite{[FST]}}]\mbox{}

  \vspace*{-4pt}
  \begin{enumerate}\addtolength{\parskip}{-9pt}
  \item There is an isomorphism of $\tSL2$ representations
    \begin{equation}
      \mV_{h,t}\tensor\Xi\approx\bigoplus_{\theta\in\oZ}\,
      \smM_{-\frac{t}{2}h,t-2;\theta}\tensor
      \mF^-_{\sqrt{\frac{t}{2}}(h+\theta)}
      \label{idspaces2}
    \end{equation}
    where on the LHS the $\tSL2$ algebra acts by
    generators~\req{invKS}, while on the RHS it acts naturally on
    $\smM_{-\frac{t}{2}h,t-2;\theta}$ as on a twisted Verma module,

  \item A singular vector exists in $\mV_{h,t}$ if and only if a
    singular vector exists in one (hence, in all) of the modules
    $\smM_{-\frac{t}{2}h,t-2;\theta}$, $\theta\in\oZ$. Whenever this
    is the case, moreover, the submodules associated with the singular
    vectors, in their own turn, satisfy an equation of the same type
    as~\req{idspaces2}.

  \end{enumerate}
\end{thm}

\vspace*{-4pt}

\noindent
The statement regarding singular vectors appeared, in a rudimentary
form, in\refnote{\cite{[S-sing]}}.  The theorem means that, as regards
the existence and the structure of submodules, the topological $\N2$
modules are {\it equivalent\/} to $\tSL2$ Verma modules: \ twisted
topological Verma submodules appear simultaneously with $\tSL2$ Verma
submodules\footnote{in particular, the $\Xi$ and $\mF^-_{\ldots}$
  modules in~\req{idspaces2} are truly `auxiliary', since nothing can
  happen there that would violate the correspondence between
  submodules in topological $\N2$ and $\tSL2$ Verma modules.}, as
\begin{equation}
  \label{compare}
  \unitlength=.6pt
\unitlength=.5pt
\begin{picture}(600,160)
  \put(0,0){
    \bezier{400}(0,0)(50,120)(118,166)
    \bezier{400}(118,166)(190,120)(240,0)
    \bezier{200}(20,0)(40,50)(70,85)
    \bezier{300}(70,85)(170,140)(230,0)
    \put(115,95){$\times$}
    \put(240,100){$\Longleftrightarrow$}
    \put(300,160){\line(1,0){180}}
    \put(480,160){\line(1,-1){120}}
    \put(310,100){\line(1,0){100}}
    \put(410,100){\line(1,-1){100}}
    \put(480,17){$\times$}
    }
\end{picture}
\end{equation}

A common feature of $\tSL2$ and topological $\N2$ Verma modules is
that all of them are generated from a state that satisfies stronger
annihilation conditions than the other states in the extremal
diagram\footnote{One {\it can\/} generate the same submodule from the
  state marked with a $\times$, but there are hardly any reasons to do
  so in the $\tSL2$ case. The point of\,\refnote{\cite{[ST4]}} is that
  doing so in the $\N2$ case is equally inconvenient.}.  What {\it
  is\/} somewhat unusual about this correspondence, though, is the
fact that on the $\tSL2$ side such a `cusp' state\,\footnote{Assigning
  grade $-n$ to $\cQ_n$ and grade $n$ to $\cG_n$, we see that every
  two adjacent arrows in diagram~\req{topdiag} represent the operators
  whose grades differ by~1, except at the cusp, where they differ
  by~2.} satisfies the same annihilation conditions as the \hw{} state
of the module, whereas on the $\N2$ side it satisfies {\it twisted\/}
topological \hw{} conditions.

Thus, to the well-known $\tSL2$ singular vectors $\ket{{\rm
    MFF}(r,s,k)}^\pm$, $r,s\in\oN$, given by the construction
of\,\refnote{\cite{[MFF]}}, there correspond the so-called {\it
  topological singular vectors\/}\refnote{\cite{[ST2],[ST3]}}
$\ket{E(r,s,t)}^\pm$, $t=k+2$, which satisfy twisted topological \hw{}
conditions \req{annihiltop} with $\theta=\mp r$ respectively:
\begin{equation}
  \cQ_{\geq\pm r}\ket{E(r,s,t)}^\pm=
  \cG_{\geq\mp r}\ket{E(r,s,t)}^\pm=
  \cL_{\geq1}\ket{E(r,s,t)}^\pm=
  \cH_{\geq1}\ket{E(r,s,t)}^\pm=0\,.
  \label{twistedtophw}
\end{equation}
As we see from the twist, the submodule generated from
$\ket{E^\pm(r,s,t)}\in\mV_{h,t}$ is the twisted
topological Verma module $\smV_{h\pm r\frac{2}{t},t;\mp r}$.
Equivalently, one may choose to describe the
positions of $\ket{E^\pm(r,s,t)}\in\mV_{h,t}$ in the
(charge, level) lattice by using the {\it eigenvalues\/} of $\cH_0$
and $\cL_0$:
\begin{equation}
  \cH_0\,\ket{E^\pm(r,s,t)}
  = h_0^\pm\,\ket{E^\pm(r,s,t)}\,,\qquad
  \cL_0\,\ket{E^\pm(r,s,t)}
  =\ell_0^\pm\,\ket{E^\pm(r,s,t)}\,,
\end{equation}
then
\begin{equation}
  h_0^\pm=h_0 \pm r\,,\qquad
  \ell_0^\pm = \ell_0 + \half r(r+2s-1)
\end{equation}
where $h_0$ and $\ell_0$ are the {\it eigenvalues\/} of $\cH_0$ and
$\cL_0$, respectively, on the \hw{} vector of the topological Verma
module~\,$\mV_{h,t}$.

The topological singular vectors occur in the topological Verma module
$\mV_{h,t}$ whenever there exist $r,s\in\oN$ such that the $h$
parameter can be represented as $h=\hminus(r,s,t)$ or
$h=\hplus(r,s,t)$, where
\begin{equation}
  \htop^-(r,s,t)=\frac{r+1}{t}-s\,,\qquad
  \htop^+(r,s,t)=-\frac{r-1}{t}+s-1
  \label{htop}
\end{equation}
The explicit construction for $\N2$ singular vectors can be found
in\refnote{\cite{[ST2],[ST3]}}.

\medskip

The idea regarding the correspondence between submodules can be
developed in the direction of category theory. Very roughly, a
category is a collection of {\it objects}, some of which may be
related by {\it morphisms}. Taking the objects to be all the (twisted)
topological $\N2$ Verma modules, the morphisms would have to be the
usual $\N2$-homomorphisms.  However, two {\it Verma\/} modules are
related by a morphism only if one of the modules can be embedded into
the other\,\footnote{Strictly speaking, this is so only for the
  ``true'' Verma modules, i.e., for the usual Verma modules on the
  $\tSL2$ side and the topological Verma modules on the $\N2$ side.
  These are precisely the modules we are interested in at the moment,
  however even for the massive $\N2$ modules and the corresponding
  relaxed $\tSl2$ modules that we introduce in what follows, the
  mappings with nontrivial kernels also agree in both categories.}. We
have just seen that such embeddings --- i.e., the occurrence of {\it
  sub\/}modules --- are `synchronized' between the topological $\N2$
Verma modules and the $\tSL2$ Verma modules. In fact, there also
exists a {\it functor\/} acting in the inverse direction, and one
eventually concludes that the category $\TOP$ of topological $\N2$
Verma modules is equivalent to the category $\VER$ of $\tSL2$ Verma
modules. To be more precise, the appearance of the twist (the spectral
flow transform) results in that this equivalence takes place only
after one effectively factorizes over the spectral flows on either
$\N2$ and $\tSl2$ sides, see\,\refnote{\cite{[FST]}} for a rigorous
statement. Anyway, an immediate consequence of this
equivalence is that\strut\\
\centerline{\framebox{\it Embedding diagrams of Verma modules are
    identical on the $\N2$ and $\tSL2$ sides,}} where we are so far
restricted to {\it topological\/} Verma modules on the $\N2$ side.
Since the $\tSL2$ embedding diagrams are
well-known\,\refnote{\cite{[RCW],[Mal]}}, this spares us the job of
deriving them in a less friendly environment of the $\N2$ algebra.


As another consequence of the equivalence theorem, the results
of\,\refnote{\cite{[MFF],[RCW],[Mal]}} reformulate as
follows:\strut\nopagebreak

\centerline{\framebox{\parbox{.96\textwidth}{\it a maximal submodule
      of a topological Verma module is either a twisted topological
      Verma module or a sum of two twisted topological Verma
      modules.}}}

\strut Since every twisted topological Verma submodule is generated
from a topological singular vector, this can be reformulated as the
statement that {\it all singular vectors in topological Verma modules
  are the topological singular vectors} (the submodules being {\it
  freely\/} generated from these vectors).  That the top-level
representative of the extremal diagram of the $\N2$ submodule
in~\req{compare} satisfies only {\it massive}, rather that
topological, \hw{} conditions does not, of course, change the fact
that the submodule is {\it not\/} a massive one: the submodule is not
{\it freely\/} generated from the top-level extremal state, as we saw
in the criterion following~\req{terminate}.  The attempts that have
been made in the literature to find massive Verma modules inside the
topological ones are erroneous.

\vspace*{-4pt}

\subsection{Where do the massive $\N2$ modules go?}

\vspace*{-4pt}

Having seen that the topological $\N2$ Verma modules are in a `good'
correspondence with $\tSL2$ Verma modules, we recall
from~\req{pic:compare} that this involves only a `small' part of $\N2$
Verma modules, whereas the massive $\N2$ modules (the `wide' ones)
seem to have nowhere to go in the $\tSL2$ picture, since all of the
capacities of the $\tSl2$ Verma modules are already used up to
maintain the correspondence with the topological (the `narrow') $\N2$
Verma modules.

\vspace*{-6pt}

\section{RELAXED $\tSL2$ VERMA MODULES}

\vspace*{-6pt}

Solving the above problem requires introducing a new class of $\tSL2$
modules. These have a characteristic property that their extremal
diagrams have {\it no `cusps'} (no angles), which will be crucial for
relating them to the massive $\N2$ Verma modules (whose extremal
diagrams have no cusps either). The recipe is to {\it relax\/} the
annihilation conditions imposed on the \hw{} vector\footnote{Yet the
  crossing out operation looked nicer in my transparencies.}:
$$
\unitlength=1pt
\begin{picture}(440,10)
  \put(-20,0){
    \put(200,0){$J^+_0\,\ketSL{j,k}~{}={}~0$}
    \linethickness{2pt}
    \bezier{50}(200,15)(250,10)(290,0)
    \bezier{50}(210,-5)(250,10)(300,15)
    }
\end{picture}
$$

For $\theta\in\oZ$, the {\it twisted relaxed Verma module\/}
$\smR_{j,\Lambda,k;\theta}$ is defined as follows. One takes the state
$\ketSL{j,\Lambda,k;\theta}$ to satisfy annihilation conditions
\begin{equation}
  J^+_{\geq\theta+1}\,\ketSL{j,\Lambda,k;\theta}=J^0_{\geq1}\,
  \ketSL{j,\Lambda,k;\theta}=
  J^-_{\geq-\theta+1}\,\ketSL{j,\Lambda,k;\theta}=0\,.
  \label{floorhw}
\end{equation}
The module is generated from $\ketSL{j,\Lambda,k;\theta}$ by a free
action of the operators $J^+_{\leq\theta-1}$, $J^-_{\leq-\theta-1}$,
and $J^0_{\leq-1}$, and by the action of operators $J^+_{\theta}$ and
$J^-_{-\theta}$ subject to the constraint
\begin{equation}\label{find-Lambda}
  J^-_{-\theta}J^+_\theta\,
  \ketSL{j,\Lambda,k;\theta}=\Lambda\,\ketSL{j,\Lambda,k;\theta}\,.
\end{equation}
In addition, the $j$ parameter is chosen such that
\begin{equation}
  \left(J^0_0+\frac{k}{2}\theta\right)\,
  \ketSL{j,\Lambda,k;\theta}=j\,\ketSL{j,\Lambda,k;\theta}
\end{equation}
(the normal ordering was chosen in~\req{find-Lambda} in order to
facilitate the evaluation of the affine Sugawara dimension of the
state).

Then, we can act on the \hw{} vector $\ketSL{j,\Lambda,k;\theta}$ with
both $J^+_0$ and $J^-_0$, thereby producing new states
\begin{equation}
  \ketSL{j,\Lambda,k;\theta|n}=\left\{\kern-4pt\new
    \begin{array}{ll}
    (J^-_{-\theta})^{-n}\,\ketSL{j,\Lambda,k;\theta}\,,&n<0\,,\\
    (J^+_\theta)^{n}\,\ketSL{j,\Lambda,k;\theta}\,,&n>0\,,
  \end{array}
\right.
  \label{theother}
\end{equation}
with $\ketSL{j,\Lambda,k;\theta|0}=\ketSL{j,\Lambda,k;\theta}$.  As a
result, the extremal diagram opens up to the straight angle; in the
case of $\theta=0$ it thus becomes
\begin{equation}
  \unitlength=1pt
  \begin{picture}(250,20)
    \put(-35,2){\Huge $\ldots$}
    \put(0,0){$\bullet$}
    \put(15,5){${}^{J^-_0}$}
    \put(28,3){\vector(-1,0){22}}
    \put(30,0){$\bullet$}
    \put(45,5){${}^{J^-_0}$}
    \put(58,3){\vector(-1,0){22}}
    \put(60,0){$\bullet$}
    \put(75,5){${}^{J^-_0}$}
    \put(88,3){\vector(-1,0){22}}
%
    \put(90,0){$\star$}
    \put(100,5){${}^{J^+_0}$}
    \put(97,3){\vector(1,0){22}}
    \put(120,0){$\bullet$}
    \put(130,5){${}^{J^+_0}$}
    \put(127,3){\vector(1,0){22}}
    \put(150,0){$\bullet$}
    \put(160,5){${}^{J^+_0}$}
    \put(157,3){\vector(1,0){22}}
    \put(180,0){$\bullet$}
    \put(193,2){\Huge $\ldots$}
  \end{picture}
  \label{floor}
\end{equation}
where all of the other states from the module correspond to points
below the line.  The $\star$ state is the above
$\ketSL{j,\Lambda,k;\theta=0}$.  We also define
$\ketSL{j,\Lambda,k|n}=\ketSL{j,\Lambda,k;0|n}$, then the norms of
these extremal states are given by
\begin{equation}
  \left\|\ketSL{j, \Lambda, k| n}\right\|^2 =
  \left\{\kern-4pt
    \begin{array}{ll}
      \prod_{i= 0}^{-n-1}(\Lambda + 2(i+1) j - i(i+1))\,,&n\leq-1\,,\\
      \prod_{i= 0}^{n-1}(\Lambda - 2 i j - i(i+1))\,, & n\geq1\,.
    \end{array}
  \right.
  \label{norms}
\end{equation}
Thus, as we move either right or left along the extremal diagram, the
norm becomes negative eventually. The negative-norm states can be
factored away if it happens that the norm of one of the extremal
states is exactly zero. This is the case whenever
$\Lambda=\Lambdach(p,j)\equiv p(p+1)+2pj$, $p\in\oZ$; then the factors
in \req{norms} become $(1 + i + p) (2 j + p - i)$ and $(p - i) (1 + i
+ 2 j + p)$ respectively. The corresponding zero-norm
state
\begin{equation}
  \ketSL{C(p,j,k)}=\left\{\kern-4pt\new\begin{array}{ll}
      (J^-_0)^{-p}\,\ketSL{j,\Lambdach(p,j),k}\,,&p\leq-1\,,\\
      (J^+_0)^{p+1}\,\ketSL{j,\Lambdach(p,j),k}\,,&p\geq0\,,
    \end{array}\right.
  \label{chargedsl2}
\end{equation}
then satisfies the Verma \hw{} conditions for $p\leq-1$ and the
twisted Verma \hw{} conditions with the twist parameter $\theta=1$ for
$p\geq1$. Thus, it is a {\it singular vector}, which can be quotiened
away along with a tail of negative-norm states.  For historical
reasons\,\refnote{\cite{[BFK],[FST]}}, states~\req{chargedsl2} are
called {\it charged\/} singular vectors---they are an $\tSL2$
counterpart of the $\N2$ singular vectors (shown schematically in the
second diagram in~\req{pic:compare}) that are called charged
since\,\refnote{\cite{[BFK]}}.

Theorem~1 is now extended to

\vspace*{-4pt}

\begin{thm}[\refnote{\cite{[FST]}}]\mbox{}

  \vspace*{-4pt}
  \begin{enumerate}\addtolength{\parskip}{-9pt}
  \item There is an isomorphism of $\tSL2$ representations
    \begin{equation}
      \mU_{h,\ell,t}\tensor\Xi\approx\bigoplus_{\theta\in\oZ}\,
      \smR_{-\frac{t}{2}h,t\ell,t-2;\theta}\tensor
      \mF^-_{\sqrt{\frac{t}{2}}(h+\theta)}
      \label{relaxedidspaces2}
    \end{equation}
    where on the LHS the $\tSL2$ algebra acts by
    generators~\req{invKS}, while on the RHS it acts naturally on
    $\smR_{-\frac{t}{2}h,t\ell,t-2;\theta}$ as on a twisted relaxed
    Verma module.

  \item A singular vectors exists in $\mU_{h,\ell,t}$ if and only if a
    singular vector exists in one (hence, in all) of the modules
    $\smR_{-\frac{t}{2}h,t\ell,t-2;\theta}$.  Whenever this is the
    case, moreover, the respective submodules, in their own turn,
    satisfy an equation of the same type as~\req{relaxedidspaces2} if
    these are massive/relaxed submodules, and Eq.~\req{idspaces2} if
    these are topological/usual-Verma submodules.
  \end{enumerate}
\end{thm}

\vspace*{-4pt}

\noindent
As follows from the notations, the parameters of the twisted relaxed
Verma module on the RHS are $j=-\frac{t}{2}h$, $\Lambda=t\ell$, and
$k=t-2$. \ The simultaneous appearance of the {\sf massive/relaxed}
and topological/usual-Verma submodules can be illustrated as follows:
$$
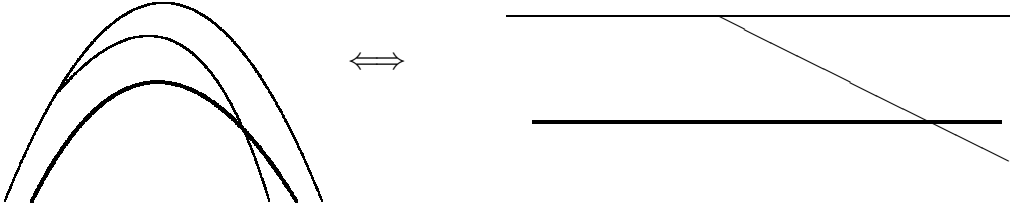
\unitlength=1pt
\begin{picture}(440,80)
  \put(20,0){
    \bezier{800}(0,0)(60,150)(120,0)
    {\linethickness{.9pt}
      \bezier{400}(10,0)(55,90)(110,0)
      }
    \bezier{500}(19,40)(70,100)(100,0)
    \put(130,50){$\Longleftrightarrow$}
    \put(190,70){\line(1,0){190}}
    \put(270,70){\line(2,-1){110}}
    {\linethickness{1.1pt}
      \put(200,30){\line(1,0){177}}
      }
    }
\end{picture}
$$
As a consequence,\strut\\
\centerline{\framebox{ \parbox[t]{.95\textwidth}{\it embedding
      diagrams of massive $\N2$ Verma modules are isomorphic to the
      embedding diagrams of relaxed $\tSL2$ Verma modules.}}}

\noindent
The analysis of the latter is easier\,\refnote{\cite{[SSi]}} because
the affine-Lie algebra representation theory is available then and
certain subdiagrams in the relaxed embedding diagrams are literally
the standard $\tSL2$ embedding diagrams\,\refnote{\cite{[RCW],[Mal]}}.
Moreover, even though the relaxed $\tSL2$ Verma modules are not a
`classical' object in the representation theory of affine Lie
algebras, the problem of enumerating submodules of relaxed Verma
modules can be reduced, to a large extent, to analysing the
standard~$\tSL2$-embedding diagrams: for a given relaxed Verma module
\,$\mR$, one can find an {\it auxiliary\/} usual-Verma module \,$\mM$
whose submodules are in a $1:1$ (or, in some degenerate cases,
essentially in a $2:1$) correspondence with the relaxed Verma
submodules in \,$\mR$, see\,\refnote{\cite{[SSi]}} for the
classification and the detailed account\,\footnote{Thus, the $\N2$
  extremal diagrams known in the
  literature\refnote{\cite{[Kir],[Dob],[M]}} need being corrected
  already for the sole reason that they do not distinguish between
  topological and massive Verma modules. That the different types of
  $\N2$ Verma-like modules were not recognized, complicates the
  analysis of degenerations of these modules
  in\,\refnote{\cite{[D]}}.}. \ Recall that for Verma modules over the
affine Lie algebras, the structure of the embedding diagrams is
governed by the affine Weyl group; for the $\N2$ algebra, we face the
problem that it is not affine and, thus, no standard construction of
an `affine' Weyl group applies. Rather, it is the
known\,\refnote{\cite{[SSi]}} $\N2$/relaxed-$\tSL2$ embedding diagrams
that should suggest the appropriate \hbox{representation of the Weyl
  group.}\nopagebreak

Let me also note that the $\N2$/relaxed-$\tSL2$ {\it embedding\/}
diagrams are made up of {\sc embeddings}, i.e., of mappings with
trivial kernels.  On the $\N2$ side, this matter appears to have
caused some confusion in the literature, because the existence of
fermions was believed to lead to the vanishing of certain would-be
embeddings. In the $\tSl2$ terms, however, this problem is obviously
absent, hence it is but an artefact on the $\N2$ side as well. The
vanishing of some compositions of the `embeddings' observed previously
is nothing but the manifestation of two facts: (i)~the criterion that
states \req{terminate} vanish for $N\gg1$ or for $M\gg1$ once
$\ket{X}$ is inside a (twisted) topological Verma module, and (ii)~the
fact that every submodule generated from a charged singular vector is
necessarily a twisted topological Verma module (similarly, and more
transparent, on the $\tSL2$ side, where the charged singular vectors
generate the usual (i.e., {\it not\/} relaxed) Verma modules, which
obviously `defermionizes' the whole picture).

Another construction which, in the affine case, reflects the structure
of the Weyl group is the BGG resolution\,\refnote{\cite{[BGG]}}.
Translating the embedding diagrams into a BGG-type resolution requires
more work in the $\N2$/relaxed-$\tSl2$ case because of the two types
of submodules existing in the appropriate Verma-like modules.  Unlike
the embedding diagrams, the resolutions are constructed in terms of
modules of only one type, therefore one would have to additionally
resolve all the twisted topological Verma modules in terms of the
massive Verma modules (or, in the $\tSL2$ terms, to resolve the
usual-Verma modules in terms of a sequence of relaxed Verma modules
with linearly growing twists).  Constructing the resolution provides
one with the tool for systematically deriving the $\N2$ characters
and finding new representations for the known characters.

\vspace{-6pt}

\section{MASSIVE AND RELAXED MODULES IN THE BOSONIC STRING}

\vspace{-6pt}

According to the above equivalence Theorems, it is inessential in many
respects whether one analyses $\N2$ Verma modules or relaxed $\tSL2$
Verma modules.  In this section, we show how the above constructions
(extremal states, massive Verma modules, etc.)  arise naturally in the
bosonic string.

In the noncritical bosonic string, one has the $\N2$ algebra realized
as in\,\refnote{\cite{[GS2],[BLNW]}}. Applying the mapping described
in the previous section, one recovers the corresponding realization of
$\tSL2$ found in\refnote{\cite{[S-sing]}}.  Let us describe in more
detail the $\N2$ version of this construction.  One starts with a
matter theory represented by the \emt{} $T$ with central
charge~$13-6/t-6t$ and tensors it with the $bc$ ghosts and a (free)
Liouville scalar. The resulting $\N2$ generators read as
\begin{equation}\new
  \begin{array}{rcl}
    \cT\kern-4pt &=&\kern-4pt T - t \d\varphi \d\varphi -
    (1+t) \d^2\varphi - \d b c - 2
    b\d c\,,\quad
    \cH = 2 \d\varphi + b c\,,\quad  \cG = b\,,\\
    \cQ\kern-4pt &=&\kern-4pt -2 b \d c c - 2 t \d\varphi \d\varphi c +
    4 \d\varphi \d c + 2 T  c +
    (2 - 2t)\d^2\varphi\, c + (1 - \frac{2}{t}) \d^2c\,,
  \end{array}
  \label{realization}
\end{equation}
where the Liouville OPE is chosen in a non-canonical normalization
$\d\varphi(z)\,\d\varphi(w)=-(1/2t)\,1/(z-w)^2$. The representation
space is then constructed as follows.

\subsection{Constructing the representation}
Each matter primary $\ketm{\Delta}$ of dimension~$\Delta$ can be
dressed into $\N2$ primaries either as
\begin{equation}
  \ketstar{h, \Delta - \frac{1}{4}(2 - 2 h - t + h^2 t), t; \theta}=
  \ket{\theta}_{\rm gh}\tensor
  \ket{e^{2t(- \half - \frac{h}{2} - \frac{\theta}{t})\varphi}}\tensor
  \ketm{\Delta}
  \label{generalstate}
\end{equation}
or by replacing $h\mapsto\frac{2}{t}-h$ in this formula (which does
not change the $\cL_0$ eigenvalue). Here,
\begin{equation}
  \ket{\theta}_{\rm gh} = \left\{ \kern-6pt 
    \begin{array}{ll}
      b\,\d b\,\ldots\,\d^{-\theta-2}b\,\ketGH{0}\,,&\theta\leq-2\,,\\
      \ketGH{0}\,,&\theta=-1\,,\\
      c\,\d c\,\ldots\,\d^{\theta}c\,\ketGH{0}\,,&\theta\geq0\,,
    \end{array}\right.
\end{equation}
are the ghost vacua in different pictures\refnote{\cite{[FMS]}}.
States~\req{generalstate} satisfy twisted massive \hw{} conditions
\req{masshw}.  Further, in the tensor product of the matter Verma
module with the ghosts (and the Liouville), each of the
states~\req{generalstate} comes together with an infinite number of
extremal states obtained by tensoring the same matter primary with
ghost vacua in different pictures:
\begin{equation}
  \label{extremal}
    \ket{\alpha}_{\rm gh}\tensor
    \ket{e^{2t(- \half - \frac{h}{2} - \frac{\theta}{t})\varphi}}
    \tensor
    \ketm{\Delta}\quad{\rm and}\quad
    \ket{\alpha}_{\rm gh}\tensor
    \ket{e^{2t(-\frac{1}{t} - \half + \frac{h}{2} -
        \frac{\theta}{t})\varphi}}
    \tensor
    \ketm{\Delta}\,,
 \quad \alpha\in\oZ\,.
\end{equation}
We thus see that in this realization,\strut\\
\centerline{\framebox{\it choosing the ghost picture corresponds to
    traveling over the extremal diagram.}}

The `bosonization' \req{realization} has the following effect:\strut\
twisted {\it topological\/} \hw{} conditions \req{annihiltop} are
satisfied whenever the dimension of a state \req{generalstate}
vanishes:
\begin{equation}
  \label{effect}
  \ketstar{h, 0, t; \theta} = \kettop{h, t; \theta}
\end{equation}
We will thus call $\ketstar{h, \ell, t; \theta}$ the {\it
  pseudomassive\/} (\hw{}) states.  The {\it generalized DDK
  prescription\/} is that there be a twisted topological primary state
among the extremal states \req{extremal}.  This is a condition on how
the parameters in the tensor product of matter and Liouville are
related: the matter dimension should be
\begin{equation}
  \Delta(h, t) = \frac{1}{4}(2 - 2 h - t + h^2 t)
  \,.
\end{equation}
Then extremal states \req{extremal} become
$$
\ket{D(h,t,\theta,\alpha)}
\equiv
\ket{\alpha}_{\rm gh}\tensor
\ket{e^{2t(- \half - \frac{h}{2} - \frac{\theta}{t})\varphi}}\tensor
\ketm{\Delta(h, t)}
=
\ketstar{h + \frac{2(\theta-\alpha)}{t},
  \frac{(\theta - \alpha)(\alpha - \theta + 1 - h t)}{t},
  t;\alpha}.
\label{Dstates}
$$
As $\alpha$ runs over the integers, the $D(h,t,\theta,\alpha)$
states fill out the extremal diagram:
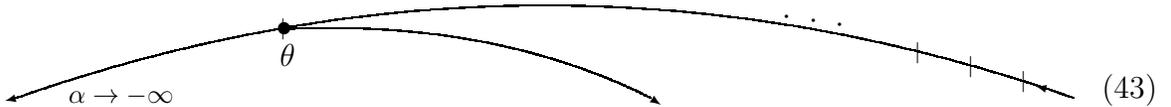
\begin{equation}
\unitlength=1pt
  \begin{picture}(400,40)
    \bezier{1000}(0,0)(200,70)(400,0)
    \put(380,4){\scriptsize $|$}
    \put(360,10){\scriptsize $|$}
    \put(340,15){\scriptsize $|$}
    \put(310,25){$\cdot$}
    \put(300,26.5){$\cdot$}
    \put(290,28){$\cdot$}
    \put(100,24){\scriptsize $|$}
    \put(98.5,23.5){$\bullet$}
    \put(100,12){$\theta$}
    \put(0,-.2){\vector(-4,-1){4}}
    \put(390,3.5){\vector(-4,1){4}}
    \put(20,-2){\footnotesize $\alpha\to-\infty$}
    \bezier{400}(100,26)(180,30)(240,0)
    \put(240,0){\vector(3,-2){4}}
  \end{picture}
  \label{alphadiagram}
\end{equation}
The states in the upper curve are freely generated from $+\infty$,
i.e., from any state where $\alpha\gg1$.  The state at $\alpha=\theta$
is the twisted topological \hw{} state~$\kettop{h,t;\theta}$, with a
twisted topological Verma submodule being generated from it. We also
have another set of extremal states
$D'(h,t,\theta,\alpha)=D(\frac{2}{t} - h,t,\theta,\alpha)$ constructed
out of the same matter $\ketm{\Delta(h, t)}$.

A straightforward analysis shows that whenever
$\ket{D(h,t,\theta,\alpha_0)}$ admits a singular vector for some
$\alpha_0\neq\theta$, each of the states $\ket{D(h,t,\theta,\alpha)}$,
$\alpha\neq\theta$, admits a massive singular vector, while
$\ket{D(h,t,\theta,\theta)}$ admits a topological singular vector.
Then the states in the $D$- and $D'$-diagrams are given by
$$
\new
\begin{array}{rcl}
  \ket{D(\hminus(r,s,t), t, \theta, \alpha)} &=&
  \ket{\alpha}_{\rm gh}\tensor
  \ket{e^{(-1 - t - r + s t - 2\theta)\varphi}}\tensor
  \ketm{\Delta_{r,s}(t)}\,,\\
  \ket{D(\hplus(r,s+1,t), t, \theta, \alpha} &=&
  \ket{\alpha}_{\rm gh}\tensor
  \ket{e^{(-1 - t + r - s t - 2\theta)\varphi}}\tensor
  \ketm{\Delta_{r,s}(t)}\,,
\end{array}\quad r,\;s\geq1\,,
$$
with $\htop^\pm(r,s,t)$ defined in~\req{htop}.  It follows,
moreover, that, using the above realization,\strut\\
\centerline{\framebox{\it\kern-3pt topological $\N2$ singular vectors
    evaluate in terms of Virasoro representation states}}\\
as follows\strut:
\begin{equation}\kern-20pt\new
  \begin{array}{crcccl}
    \begin{array}{r}E^+(r,1,t)\\
      r\geq1\end{array}&\begin{array}{r}E^+(r,s+1,t)\\
      r,s\geq1
    \end{array}\kern-30pt
    &{}&{}&{}&\kern-30pt
    \begin{array}{l}E^-(r,s,t)\\
      r,s\geq1
    \end{array}\kern-40pt\\
    \Bigm\downarrow&{}&
    \mbox{\Large$\searrow$}&{}&\mbox{\Large$\swarrow$}\\
    \mbox{Virasoro \hw{} state}&{}&{}&
    \kern-32pt \mbox{Virasoro singular vector}\,(r,s|\Delta_{r,s})
    \kern-32pt&{}
  \end{array}
  \label{reduction}
\end{equation}
In view of the equivalence theorem, the same reduction holds for
singular vectors in the usual-Verma $\tSL2$ modules, as known since
long ago\,\refnote{\cite{[GP]}}.  The fact that the $\ket{E(r,1,t)}^+$
topological singular vectors cannot be constructed out of the matter
ones shows up in a different guise as the existence of a screening
operator in the realization~\req{realization} (or, equivalently, a
similar realization\refnote{\cite{[S-sing]}} of the $\tSL2$ algebra).

\subsection{Mappings by the screening current}
Whenever one uses an operator construction (`bosonization'), leading
to some additional effects (e.g., vanishings) in representations
(Eq.~\req{effect} in our case), one should expect the appearance of a
screening current.  We, indeed, have a {\it fermionic screening
  current\/} of the form\refnote{\cite{[scr],[q]}}
\begin{equation}
  F=b\,e^{t\varphi}\,\Psi_{12}\,,
\end{equation}
where $\Psi_{1,2}$ is the `$12$' operator in the matter (Virasoro)
 sector. The  $\Psi_{1,2}$ operator has
two components that can be distinguished by picking out the following
terms from the fusion relations:
\begin{equation}
  \Psi_{12}^\pm*\ketm{\Delta(h,t)}\sim\ketm{\Delta(h\pm1,t)}\,.
  \label{PsionVerma}
\end{equation}
Here, the \hw{} states may be understood as those in {\it Verma\/}
modules over the Virasoro algebra.  We now can construct the following
action of the screening $F$ on the pseudomassive modules over the
$\N2$ algebra (we omit the integral which makes the screening {\it
  charge\/} out of the current):
\begin{equation}
  F^\pm : \ket{D(h,t,\theta,\alpha)}\mapsto
  \left\{\kern-4pt
    \new\begin{array}{ll}
      \ket{D(h\pm1,t,\theta,\alpha-1)}\,,&
      \alpha\geq\theta + 1\,,\\
      0\,,&{\rm otherwise}
    \end{array}
  \right.
\end{equation}
(which of the extremal states do, and which do not, vanish under the
action of the screening, follows from a simple analysis of operator
products). Then, for $\alpha\geq\theta+1$,
\begin{equation}\new
  \begin{array}{rcrcl}
    F^-\kern-5pt&:&\kern-5pt
    \ket{D(\hminus(r,s,t),t,\theta,\alpha)}&\mapsto&
    \ket{D(\hminus(r,s+1,t),t,\theta,\alpha-1)}\,,\\
    F^+\kern-5pt&:&\kern-5pt
    \ket{D(\hplus(r,s+1,t),t,\theta,\alpha)}&\mapsto&
    \ket{D(\hplus(r,s,t),t,\theta,\alpha-1)}\,.
  \end{array}
  \label{2mappings}
\end{equation}

Next, we observe that the identities
\begin{equation}
  \ket{D(\hminus(r,s,t), t, \theta_1, \alpha)}=
  \ket{D(\hplus(r,s+1,t), t, \theta_2, \alpha)}\,,\quad\forall\alpha\,,
  \label{key}
\end{equation}
hold if and only if either $s=0$, $r + \theta_1 = \theta_2$, or $t={r
  + \theta_1 - \theta_2\over s}$. For generic (non-rational) $t$, we
can use the $s=0$ case in order to connect the two series of mappings
\req{2mappings} together.  Omitting the $t$ parameter, we label the
extremal diagrams spanned out by the
$\ket{D(\htop^\pm(r,s,t),t,\theta,\alpha)}$ states by the
corresponding $\htop^\pm(r,s)$ and the value of $\theta$ that gives
the position of the topological \hw{} vector. We then have the
following mappings of modules with the extremal
diagrams~\req{alphadiagram}
\begin{equation}
  \begin{picture}(400,30)
    \put(0,-20){
      \put(40,0){
        \put(-25,35){${}^{F^+}$}
        \put(-45,30){$\ldots\longrightarrow$}
        \bezier{100}(0,30)(20,40)(40,30)
        \put(30,31){\scriptsize $\bullet$}
        \put(0,44){${}_{\hplus(r, 2), r + \theta}$}
        }
      \put(125,0){
        \put(-25,35){${}^{F^+}$}
        \put(-30,30){$\longrightarrow$}
        \bezier{100}(0,30)(20,40)(40,30)
        \put(10,32){\scriptsize $\bullet$}
        \put(30,31){\scriptsize $\bullet$}
        \put(0,44){${}_{\hplus(r, 1), r + \theta}$}
        \put(5,20){${}_{\hminus(r, 0), \theta}$}
        }
      \put(210,0){
        \put(-25,35){${}^{F^-}$}
        \put(-30,30){$\longrightarrow$}
        \bezier{100}(0,30)(20,40)(40,30)
        \put(10,32){\scriptsize $\bullet$}
        \put(5,20){${}_{\hminus(r, 1), \theta}$}
        }
      \put(295,0){
        \put(-25,35){${}^{F^-}$}
        \put(-30,30){$\longrightarrow$}
        \bezier{100}(0,30)(20,40)(40,30)
        \put(10,32){\scriptsize $\bullet$}
        \put(5,20){${}_{\hminus(r, 2), \theta}$}
        }
      \put(380,0){
        \put(-25,35){${}^{F^-}$}
        \put(-30,30){$\longrightarrow\ldots$}
        }
      }
  \end{picture}
  \label{glued}
\end{equation}

This sequence applies to those modules where the Virasoro part is
taken to be {\it Verma\/} modules. We now investigate whether it is
possible to go over from \req{glued} to a similar sequence for {\it
  quotient\/} modules.  According to~\req{reduction}, the massive
$\N2$ singular vectors would be factored away in all of the terms
starting with and after $\hminus(r, 1)$ as soon as the Virasoro
singular vectors are factored away. This allows one to define the
$F^-$ mappings between the {\it irreducible\/} representations.  The
same is true for the modules before and including $\hplus(r,2)$.  In
the middle term $\hminus(r,1)$, however, {\it there is no submodule\/}
to factor over in the corresponding Virasoro Verma module. Yet, taking
the composition $F^-\circ F^+$ allows us to conclude
that\\
\centerline{\framebox{\it for $t\in\oC\setminus\oQ$, there exists an
    exact sequence}}
\begin{equation}\new
  \begin{array}{l}
    \kern-20pt\ldots\,\stackrel{F^+}{\longrightarrow}
    \ket{\bar{D(\hplus(r,3,t),t,\theta+r)}}
    \stackrel{F^+}{\longrightarrow}
    \ket{\bar{D(\hplus(r,2,t),t,\theta+r)}}
    \stackrel{F^-\circ F^+}{\longrightarrow}{}\\
    \kern120pt
    \ket{\bar{D(\hminus(r,1,t),t,\theta)}}
    \stackrel{F^-}{\longrightarrow}
    \ket{\bar{D(\hminus(r,2,t),t,\theta)}}
    \stackrel{F^-}{\longrightarrow}\,\ldots
  \end{array}
\end{equation}
where the bars indicate that the Virasoro singular vectors are
declared to vanish (i.e., {\it irreducible\/} representations are
taken in the matter sector in~\req{generalstate}, \req{extremal}, and
similar formulae). The exact sequences of this sort (actually, those
involving Verma modules) make up a part of the BGG resolution for the
irreducible $\N2$ representations in terms of the massive Verma
modules. As I have mentioned, the BGG resolution allows one to find the
characters; while for the affine Lie algebras the thus found character
formulae reproduce the Weyl--Ka\v c formula, the importance of the
present approach consists in that, with the appropriate representation
of the Weyl group not known, the BGG resolution has to be constructed
directly from the corresponding embedding diagrams.

Another useful observation is that the above exact sequence is
parallel to an exact sequence between representations of the quantum
group $s\ell(2|1)_q$ (for $q$ not a root of unity in accordance with
the above choice of generic $t$), which is not a
coincidence\refnote{\cite{[q]}}. The $s\ell(2|1)_q$ embeddings are
also performed by `fermionic' singular vectors, with a due analogue of
the $F^-\circ F^+$ composition in the center.  In fact, generating the
quantum $s\ell(2|1)_q$ symmetry involves other screenings in addition
to $F$, however these do not have to be explicitly introduced in the
present approach, where we do not bosonize the matter (\emt~$T$)
through {\it free\/} fields.  The $s\ell(2|1)$ quantum group has long
been suspected to govern the complete $\tSL2$ fusion
rules\,\refnote{\cite{[FM]}}, however the presently observed symmetry
is only~$osp(1|2)_q$. It may be expected that the model with a
$s\ell(2|1)_q$-symmetric fusion may be constructed by developing the
above observations.

\subsubsection{Acknowledgements.} It is a pleasure to thank the
organizers of the conferences, first of all P.H.~Damgaard and
J.~Jurkiewicz, and A.~Belavin. I am grateful to B.~Feigin and
I.~Tipunin for a fruitful collaboration, and to M.~Halpern,
F.~Malikov, I.~Shchepochkina, V.~Sirota, and V.~Tolstoy for useful
discussions. This work was supported in part by RFFI Grant~96-02-16117.

\vfill

\end{document}